\documentclass[reqno,10pt,a4paper]{amsart}
\usepackage[margin=2.7cm]{geometry}
\usepackage{cite}
\usepackage{newpxtext}
\usepackage{cabin}% sans serif
\usepackage[varqu,varl]{inconsolata}% sans serif typewriter
\usepackage[bigdelims,vvarbb]{newpxmath}% bb from STIX
\usepackage[cal=cm,bb=esstix,bbscaled=1.05]{mathalfa}% mathcal
\usepackage{graphicx}
\usepackage{colortbl}
\usepackage{booktabs}
\usepackage[up]{caption}
\usepackage[labelformat=simple]{subcaption}
\usepackage{xcolor}
\usepackage{hyperref}
\hypersetup{
    colorlinks,
    linkcolor={red!50!black},
    citecolor={blue!50!black},
    urlcolor={blue!80!black}
}
\AtBeginDocument{%
\def\MR#1{}
}
\newcommand{\orcid}[1]{\,\resizebox{8px}{!}{\href{https://orcid.org/#1}{\includegraphics{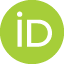}}}}
\graphicspath{{images/}}
\newcommand{\ZZ}{\mathbb{Z}}
\newcommand{\PP}{\mathbb{P}}
\newcommand{\wPP}{w\PP}
\newcommand{\CC}{\mathbb{C}}
\newcommand{\QQ}{\mathbb{Q}}
\newcommand{\RR}{\mathbb{R}}

\newcommand{\Ghat}{\widehat{G}}

\DeclareMathOperator{\Proj}{Proj}
\DeclareMathOperator{\codim}{codim}

\DeclareMathOperator{\Hom}{Hom}
\newcommand{\NQ}{N_\QQ}
\newcommand{\MQ}{M_\QQ}
\newcommand{\QQFano}{$\QQ$\nobreakdash-Fano}

\DeclareMathOperator{\Gr}{Gr}
\DeclareMathOperator{\Ehr}{Ehr}
\newcommand{\wGr}{w\!\Gr}

\newcommand{\grdb}[1]{\href{http://grdb.co.uk/search/fano3?id_cmp=in\&id=#1}{#1}}
\theoremstyle{plain}

\theoremstyle{definition}

\theoremstyle{remark}

\numberwithin{equation}{section}
\begin{document}
%-------------------------------------------------------------------------------
\author[G.\,Brown]{Gavin Brown\orcid{0000-0002-4087-5624}}
\address{Mathematics Institute\\Zeeman Building\\University of Warwick\\Coventry\\CV4 7AL\\UK}
\email{G.Brown@warwick.ac.uk}% Brown
\author[T.\,Coates]{Tom Coates\orcid{0000-0003-0779-9735}}
\address{Department of Mathematics\\Imperial College London\\180 Queen's Gate\\London\\SW7 2AZ\\UK}
\email{t.coates@imperial.ac.uk}% Coates
\author[A.\,Corti]{Alessio Corti\orcid{0000-0002-9009-0403}}
\address{Department of Mathematics\\Imperial College London\\180 Queen's Gate\\London\\SW7 2AZ\\UK}
\email{a.corti@imperial.ac.uk}% Corti
\author[T.\,Ducat]{Tom Ducat\orcid{0000-0002-6153-6293}}
\address{Department of Mathematical Sciences\\Durham University\\Stockton Road\\Durham\\DH1 3LE\\UK}
\email{thomas.ducat@durham.ac.uk} % Ducat
\author[L.\,Heuberger]{Liana Heuberger\orcid{0000-0002-8038-9620}}
\address{D\'epartement math\'ematiques\\Facult\'e des Sciences\\2 Boulevard Lavoisier\\Angers\\France}
\email{liana.heuberger@univ-angers.fr}% Heuberger
\author[A.\,M.\,Kasprzyk]{Alexander Kasprzyk\orcid{0000-0003-2340-5257}}
\address{School of Mathematical Sciences\\University of Nottingham\\Nottingham\\NG7 2RD\\UK}
\email{a.m.kasprzyk@nottingham.ac.uk}% Kasprzyk
%-------------------------------------------------------------------------------
\thanks{\emph{Funding.} TC is supported by ERC Consolidator Grant~682603 and EPSRC Programme Grant~EP/N03189X/1. AC is supported by EPSRC Programme Grant~EP/N03189X/1. LH is supported by Projet \'Etoiles Montantes GeBi de la R\'egion Pays de la Loire. AK is supported by EPSRC Fellowship~EP/N022513/1.}
%-------------------------------------------------------------------------------
%\subjclass[2010]{14J33, 52B20 (Primary); 14J45, 14N35, 13F60, 32G20 (Secondary)}
%-------------------------------------------------------------------------------
\title{Computation and Data in the Classification of Fano Varieties}
%-------------------------------------------------------------------------------
\begin{abstract}
Fano varieties are `atomic pieces' of algebraic varieties, the shapes that can be defined by polynomial equations. We describe the role of computation and database methods in the construction and classification of Fano varieties, with an emphasis on three-dimensional Fano varieties with mild singularities called \QQFano{} threefolds. The classification of \QQFano{} threefolds has been open for several decades, but there has been significant recent progress. These advances combine computational algebraic geometry and large-scale data analysis with new ideas that originated in theoretical physics.
\end{abstract}
%-------------------------------------------------------------------------------
\maketitle
%-------------------------------------------------------------------------------
\section{Introduction}
%-------------------------------------------------------------------------------
Fano varieties are an important class of geometric shapes. They are basic building blocks in algebraic geometry, both via the Minimal Model Program~\cite{Mor82,Mor88,BCHM10,HM10} and as a source of explicit constructions. For example, one can construct three-dimensional Calabi--Yau manifolds (smooth varieties with zero curvature) as anticanonical sections in smooth four-dimensional Fano varieties; these are of particular importance in constructing models of spacetime in Type~II string theory~\cite{AGM94,Pol05}. The classification of Fano varieties is a fundamental problem in geometry yet, despite almost a century of study, this classification is still far from understood. In this paper we discuss the classification of Fano varieties over $\CC$. In particular, the $n$-dimensional Fano varieties that we consider can also be thought of, at least when they are smooth, as $2n$-dimensional real manifolds with positive curvature.

There are finitely many deformation families of Fano manifolds (that is, smooth Fano varieties) in each dimension~\cite{KMM92}. There is exactly~$1$ one-dimensional Fano manifold: the Riemann sphere~$\PP^1$. In two dimensions there are~$10$ deformation families, the del~Pezzo surfaces:~$\PP^1\times\PP^1$,~$\PP^2$, and~$\PP^2$~blown-up in at most eight points~\cite{delPez87}. In dimension three the classification is due to Fano~\cite{Fan47}, Iskovskikh~\cite{Isk77,Isk78,Isk79}, and Mori--Mukai~\cite{MM81,MM03}: there are~$105$ deformation families\footnote{The three-dimensional classification can be explored online via the excellent Fanography website~\cite{Bel}.}. Very little is known about the classification of Fano manifolds in dimensions four or more.

When we allow mild singularities even less is known. There are finitely many deformation families in each dimension with singularities of bounded complexity~\cite{KMMT00,Bir21}, but there are few explicit classification results. In dimension three it is natural to focus on the \emph{\QQFano{} threefolds}: Fano varieties with~$\QQ$-factorial terminal singularities -- this is the natural setting for the Minimal Model Program. Even here there is as yet no classification, although several hundred deformation families have been described; see e.g.~\cite{CF93,San95,San96,Fle00,Tak02,Kas06,BKR12,PR16,BKQ18,Duc18,CD20}. What is known, however, are the possible Hilbert series of \QQFano{} threefolds\footnote{Strictly speaking, we mean here the possible Hilbert series of Mori--Fano threefolds only -- see~\cite[\S5]{BK22} for a precise discussion -- but as every known \QQFano{} threefold has a Hilbert series that is in this list, we will blur the distinction in the rest of the paper. In any case, this is a large and important part of the landscape of \QQFano{} threefolds.}
~\cite{ABR02,BK22}. The Hilbert series of a Fano variety~$X$ is the generating series for the dimensions of the graded pieces of the anticanonical ring
\begin{equation}\label{eq:anticanonical}
R(X, {-K}_X) = \bigoplus_{n = 0}^\infty H^0(X, {-n} K_X).
\end{equation}
Hilbert series remain invariant under~$\QQ$-Gorenstein deformation of~$X$.  The possible Hilbert series are recorded in the Graded Ring Database~\cite{grdb,fanodata}. One can think of this database as a first sketch, giving numerical invariants only, of the possible `geography' of \QQFano{} threefolds. 

\begin{figure}[htbp]
\centering
\includegraphics[width=0.7\textwidth]{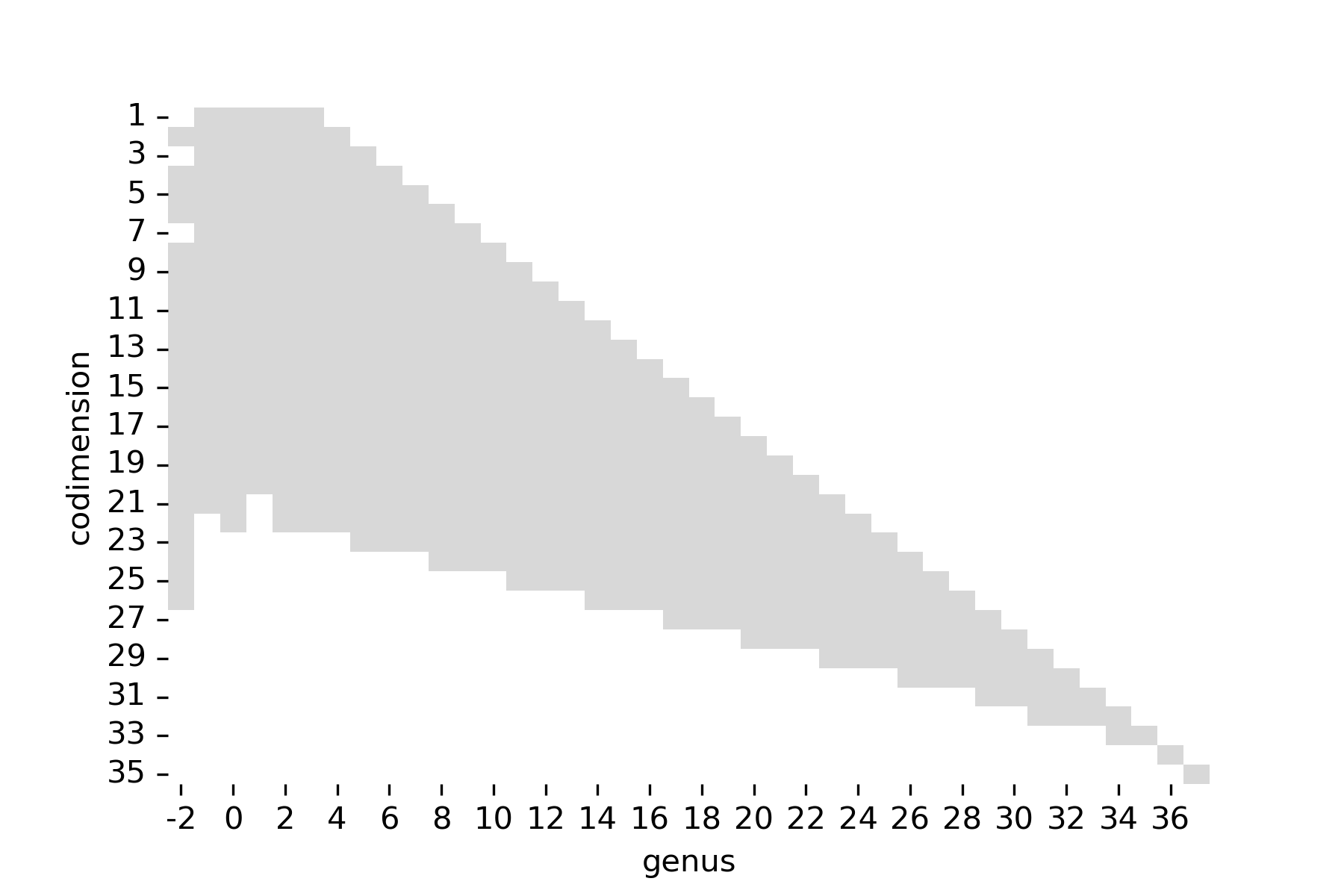}
\caption{A sketch of the landscape of \QQFano{} threefolds, from the Graded Ring Database~\cite{fanodata}.}
\label{fig:grdb}
\end{figure}

A Fano variety~$X$ can be recovered as~$\Proj$ of its anticanonical ring~\eqref{eq:anticanonical}. Choosing homogeneous generators for the anticanonical ring defines an embedding of~$X$ into a weighted projective space~$\wPP$, and we can read off various measures of the complexity of this embedding (or of the complexity of~$R(X, {-K}_X)$ as a ring) directly from the Hilbert series. In particular, the Hilbert series can be used to estimate a lower bound for the~\emph{codimension}
\[
 \codim{X} \coloneqq \dim{\wPP}-\dim{X},
\] 
that is, the number of independent equations required to define~$X$ (locally) inside~$\wPP$. The Hilbert series also determines the~\emph{genus} of the embedding
\[
 g(X) \coloneqq \dim H^0(X,{-K}_X)-2.
\]
If~$X$ is smooth then~$g(X)$ is the genus of the curve defined by intersecting~$X$ with two generic hyperplanes in~$\wPP$; for the general case, including an interpretation of the case where~$g(X)<0$, see~\cite[\S4]{ABR02}. The genus increases with the complexity of the embedding. Figure~\ref{fig:grdb} shows the range of pairs~$(g(X), \codim{X})$ for entries in the database of possible Hilbert series of \QQFano{} threefolds. See~\cite{BK22} for precise details: in particular, codimension here means estimated codimension as there.

The database~\cite{grdb,fanodata} is produced using combinatorial methods that do not guarantee existence of a \QQFano{} threefold that realises any given entry: there can be zero, one, or many deformation families of \QQFano{} threefolds with a given candidate Hilbert series. In the remainder of the paper, we discuss techniques for constructing and classifying \QQFano{} threefolds, guided by the geography in Figure~\ref{fig:grdb}.
%-------------------------------------------------------------------------------
\section{Classifications in low codimension}\label{sec:low_codim}
%-------------------------------------------------------------------------------
The first examples of Fano varieties are projective spaces~$\PP^n$. After these, projective hypersurfaces of low degree give the archetypal Fano varieties: a cubic surface, a cubic or quartic threefold, and so on. Broadening our view to orbifolds with cyclic quotient singularities greatly extends the range of hypersurfaces that we can construct. For example, a weighted homogeneous polynomial~$f$ of degree~$d$ in~$n+1$ variables with positive integer weights~$a_0,\dots,a_n$ defines a hypersurface~$(f=0)$ in weighted projective space~$\PP(a_0,\dots,a_n)$. If $d < a_0 + \cdots + a_n$ then this hypersurface is Fano. See~\cite{Mor75,Dol82,Fle00}.

In three dimensions, there are exactly~$95$ families of \QQFano{} threefolds that arise as quasismooth weighted projective hypersurfaces~\cite{Rei80,JK01,BK16b}. Here \emph{quasismoothness} ensures that all singularities of the hypersurface arise from singularities of the ambient weighted projective space~\cite{Fle00}. The $95$~hypersurfaces are all of the form~$X_d\subset\PP(1,a_1,a_2,a_3,a_4)$, where~$d=a_1+a_2+a_3+a_4$, and correspond to the~$95$ families of~K3 surfaces~$X_d\subset\PP(a_1,a_2,a_3,a_4)$. Quasismooth members $X$ of each family are birationally rigid~\cite{CPR00,CP17}, meaning if $X$ is birational to another \QQFano{} threefold or Mori fiber space $Y$ then $X$ is in fact biregular isomorphic to $Y$. One may consider this to be a strong form of non-rationality~\cite{IM71} since, in particular, they do not admit birational maps to~$\PP^3$. The~$95$ families have been studied from many points of view, and serve as a testing ground for conjectures and computations. For example, Kim--Okada--Won compute their log canonical thresholds, and deduce~$K$-stability results for~$50$ of the~$95$ families~\cite{KOW20}. In four dimensions, \QQFano{} quasismooth weighted projective hypersurfaces have also been classified~\cite{BK16}: there are~$11\,618$ families. 

Hypersurfaces~$(f=0)$ have codimension one, because they are defined (globally) by a single equation. One can construct Fano varieties with higher codimension by considering complete intersections
\[
    (f_1 = \cdots = f_c = 0)
\]
in weighted projective spaces, where the~$f_i$ are homogeneous polynomials of sufficiently low degree. Iano-Fletcher describes~$85$ families of \QQFano{} threefolds that arise as codimension-two complete intersections in weighted projective space~\cite{Fle00}; this list is in fact complete~\cite{CCC11}. In higher codimension most Fano varieties are not weighted projective complete intersections. In codimension three, for example, sections of weighted Grassmannians~$\wGr(2,5)$, defined by five equations of Pfaffian type, are typical~\cite{CR02}. Codimension four embeddings are also tractable~\cite{Alt98,BKR12} but beyond codimension four, a comprehensive understanding is restricted to special classes such as smooth Fano threefolds.
%-------------------------------------------------------------------------------
\section{Birational geometry and steps into higher codimension}\label{sec:unproj}
%-------------------------------------------------------------------------------
Whilst the low codimension studies discussed in~\S\ref{sec:low_codim} are very successful and comprehensive, they break down in higher codimension. There are several approaches to this problem. One is to find simple examples whose anticanonical embedding happens to be in high codimension: for example,~$\PP^3$ itself has~$-K_X = 4H$, for a linear hyperplane~$H\subset\PP^3$, so that the anticanonical embedding is the 4th Veronese embedding~$\PP^3\hookrightarrow\PP^{34}$ in codimension~$31$. Another approach is to find varieties that naturally live in high codimension whose sections are \QQFano{} threefolds~\cite{QS11,BKZ19}. A more sophisticated approach, the \emph{Takeuchi program}~\cite{Tak89}, analyses possible birational links between Fano threefolds with the aim of constructing complicated examples from simpler ones; cf.~\cite{Ale94,Cor95}. This was exploited by Takagi~\cite{Tak02} to classify all \QQFano{} threefolds~$X$ having only~$\ZZ/2$ hyperquotient singularities, Picard rank~$\rho_X=1$, and~$\dim H^0(X,-K_X)\ge4$; and by Prokhorov~\cite{Pro10} to analyse \QQFano{} threefolds with~$-K_X$ divisible.

We discuss a fourth approach here, applying \emph{unprojection}~\cite{PR04,Rei02,BKR12}. This is close in spirit to the Takeuchi program, and is rooted in Fano's original approach to the classification problem via projection. The idea is as follows. If~$X_4\subset\PP^4$ is a terminal Fano quartic threefold that contains a linear plane~$\PP^2\cong D\subset X_4$, then there is a birational map~$X_4\dashrightarrow Y=Y_{3,3}\subset\PP(1,1,1,1,1,2)$. In coordinates~$x_0,\dots x_4$ on~$\PP^4$, this is
\[
D = (x_0=x_1=0) \subset X_4 = (x_0A - x_1B=0) \subset \PP^4
\]
where~$A$ and~$B$ are cubic forms, and setting~$y = A/x_1 = B/x_0$ defines~$Y$ and a morphism
\begin{equation}\label{eq:unproj}
\begin{array}{r@{\ }c@{\ }l}
X_4 & \dashrightarrow &Y = (yx_1=A, yx_0=B) \subset \PP(1,1,1,1,1,2) \\
(x_0,\dots,x_4) & \mapsto &(x_0,\dots,x_4,y(x_0,\dots,x_4)).
\end{array}
\end{equation}
Here~$x_0,\ldots,x_4$ and~$y$ are the coordinates on~$\PP(1,1,1,1,1,2)$ of weights~$1,\ldots,1$ and~$2$. The map \eqref{eq:unproj} is called the \emph{unprojection} of~$D\subset X_4$, and conversely, elimination of the variable~$y$ from the equations of~$Y$ recovers~$X_4$ as a projection of~$Y$. 

Multiplying up the denominator in~\eqref{eq:unproj} of either expression for~$y$ shows that general points of~$D$ are mapped to~$[0\!:\!0\!:\!0\!:\!0\!:\!0\!:\!1]\in Y$, and so~$D$ is indeed contracted. However, the map \eqref{eq:unproj} is not a morphism. The geometry of the map depends on the singularities that~$X_4$ acquires at the points $(A=B=0)$ on~$D$. Typically these are ordinary double points, and the map consists of the~$D$-ample small resolution of  each such point followed by the contraction of~$D$, but this requires detailed analysis in any given case. If one can repeat this with two planes~$D_1$,~$D_2\subset X_4$, unprojecting them one after the other, the result would be a Fano threefold embedded in codimension three. It is easy to imagine continuing this process, but imposing successive planes becomes an increasingly difficult geometrical and combinatorial problem.

This approach works well into codimension four, with several hundred families constructed by unprojecting weighted planes~$D=\PP(1,a,b)$~\cite{BKR12}. In the example above, Brown--Kerber--Reid construct three different deformation families of \QQFano{} threefolds, using different configurations of planes. (Takagi's classification describes two of these families; the third has higher Picard rank.) Unprojection techniques have been extended into codimension~$\ge5$ by Brown and Ducat~\cite{BD22}. The real power of this technique becomes apparent when combined with ideas from mirror symmetry: see~\S\ref{sec:2pronged} below.
%-------------------------------------------------------------------------------
\section{Polytopes, fans, and toric varieties}
%-------------------------------------------------------------------------------
Toric varieties are highly symmetrical algebraic varieties. An~$n$-dimensional normal projective variety~$X$ is a~\emph{toric variety} if it contains a dense algebraic torus~$(\CC^\times)^n$ and the action of that torus on itself extends to an action of the torus on all of~$X$: see~\cite{Dan78,Ful93}. For example,  projective spaces and weighted projective spaces are toric varieties. In our context, toric varieties will play two different roles. Three-dimensional toric varieties will occur as degenerations of \QQFano{} threefolds; in general these toric varieties will be highly singular. Higher-dimensional toric varieties will occur as ambient spaces when we construct \QQFano{} threefolds as hypersurfaces or complete intersections; we have already seen this in~\S\ref{sec:low_codim}.

The study of toric varieties is underpinned by a rich `dictionary' that translates geometry into combinatorics. Let~$M\cong\ZZ^n$ denote the lattice of characters of~$T = (\CC^\times)^n$, with dual lattice $N\coloneqq\Hom(M,\ZZ)$. A toric variety~$X$ with dense torus~$T$ has a combinatorial description in terms of a fan~$\Sigma$ in~$\NQ$, and many geometric properties of~$X$ can be rephrased in terms of combinatorial statements about~$\Sigma$. For example, let~$\rho_1,\ldots,\rho_k$ be the one-dimensional cones of~$\Sigma$. Each $\rho_i$ is generated by a (unique) primitive lattice element~$v_i\in N$, and~$X$ is Fano if and only if~$v_1,\ldots,v_k$ correspond to the vertices of a convex lattice polytope in~$\NQ$ containing the origin in its interior. Such a polytope is called a~\emph{Fano polytope}. Conversely, given a Fano polytope~$P$ the fan~$\Sigma$ (and hence the toric Fano variety~$X$) can be recovered by taking the collection of cones spanned by the faces of~$P$; we call this the~\emph{spanning fan} of~$P$. Under this dictionary there is a correspondence between the dual polytope
\[
    P^*\coloneqq\{u\in\MQ\mid u(v)\geq -1\}
\] 
of~$P$ and the anticanonical divisor~${-K}_X$ of~$X$. Since~$P^*$ always contains the origin, it follows that~$\dim{H^0(X,{-K}_X)}=|P^*\cap M|$ is strictly positive.
%%-------------------------------------------------------------------------------
\section{Mirror symmetry}\label{sec:mirror_symmetry}
%-------------------------------------------------------------------------------
A new approach to the classification of Fano varieties centres around a set of ideas from string theory called~\emph{mirror symmetry}~\cite{GP90,COGP91,HV00}. Mirror symmetry gives a correspondence between Fano manifolds and certain Laurent polynomials~\cite{CCGGK13}. From this perspective, the key invariant of a Fano variety~$X$ is its~\emph{regularised quantum period}
\begin{equation}\label{eq:Ghat}
\Ghat_X(t)=\sum_{d=0}^\infty c_d t^d.
\end{equation}
This is a power series with coefficients~$c_0 = 1$,~$c_1 = 0$, and~$c_d = r_d d!$, where~$r_d$ is a certain Gromov--Witten invariant of~$X$. Intuitively speaking,~$r_d$ is the number of rational curves in~$X$ of degree~$d \geq 2$ that pass through a fixed generic point and have a certain constraint on their complex structure. Under mirror symmetry, an~$n$-dimensional Fano manifold~$X$ corresponds to a Laurent polynomial $f \in \CC[x_1^{\pm 1}, \ldots, x_n^{\pm 1}]$ if the regularised quantum period~$\Ghat_X$ agrees with the~\emph{classical period}
\begin{equation}\label{eq:pi}
\pi_f(t) = \sum_{d=0}^\infty c'_d t^d.
\end{equation}
Here~$c_d'$ is the coefficient of the unit monomial in~$f^d$.

It is expected that if a Fano manifold~$X$ corresponds under mirror symmetry to a Laurent polynomial~$f$, then there is a~$\QQ$-Gorenstein degeneration of~$X$ to a (singular) toric variety~$X_f$. This toric variety is defined by the spanning fan of the Newton polytope of~$f$. One can hope, therefore, to recover~$X$ directly from its mirror partner~$f$: by regarding~$f$ as encoding a log structure on~$X_f$ that determines the smoothing from~$X_f$ to~$X$, and then constructing~$X$ as this smoothing. 

Gromov--Witten invariants are deformation invariant, so the regularised quantum period~$\Ghat_X$ is a deformation invariant of~$X$. The classical period~$\pi_f$, on the other hand, is invariant under an equivalence relation called~\emph{mutation}~\cite{ACGK12}. It is conjectured that mirror symmetry gives a one-to-one correspondence between deformation-equivalence classes of \QQFano{} varieties\footnote{That is, Fano varieties with~$\QQ$-factorial terminal locally toric singularities.} that admit a toric degeneration, and mutation-equivalence classes of~\emph{rigid maximally mutable Laurent polynomials}~(rigid MMLPs)~\cite{CKPT21}. For the classification of two- and three-dimensional Fano manifolds from this perspective, see~\cite{CCGK16,ACCHKOPPT16,DH16,KNP17,Pri18}.

If~$X$ is a Fano projective hypersurface then the Givental/Hori--Vafa mirror construction~\cite{Giv98,HV00} gives rise to a Laurent polynomial~$f$ that corresponds to~$X$ under mirror symmetry. For example, if~$X$ is a cubic surface in~$\PP^3$ then the Givental/Hori--Vafa mirror is given by 
\begin{align*}
W = x_0 + x_1 + x_2 + x_3 && \text{where} &&&
x_0 x_1 x_2 x_3 = 1 \\
&& &&& x_1 + x_2 + x_3 = 1.
\end{align*}
Eliminating the second equation by setting 
\begin{align*}
x_1 = \frac{1}{1+x+y} &&
x_2 = \frac{x}{1+x+y} &&
x_3 = \frac{y}{1+x+y} 
\end{align*}
and solving for~$x_0$ yields the Laurent polynomial
\[
f = W - 7 = \frac{(1+x+y)^3}{xy} - 6.
\]
This is the Przyjalkowski method. It applies in greater generality to give Laurent polynomials that correspond under mirror symmetry to complete intersections in weighted projective spaces~\cite{Prz11} and toric varieties~\cite{CKP15,DH16} that satisfy mild combinatorial conditions.

The Givental/Hori--Vafa construction above can be reversed to recover a presentation of~$X$ as a toric complete intersection from the corresponding Laurent polynomial~$f$: this is~\emph{Laurent inversion}~\cite{CKP19}. This process cannot always work -- for example, only~$89$ of the~$105$ smooth Fano threefolds are known to be toric complete intersections -- but nonetheless Laurent inversion gives a powerful tool for constructing \QQFano{} varieties. In outline: one starts with a Fano polytope~$P$, searches for rigid MMLPs~$f$ with Newton polytope~$P$, and then applies Laurent inversion to construct a corresponding Fano variety. Heuberger has used this approach to construct many new \QQFano{} threefolds~\cite{Heu22}.
%-------------------------------------------------------------------------------
\section{Fano polytopes}\label{sec:high_codim}
%-------------------------------------------------------------------------------
The classification of Fano polytopes is an attractive combinatorial problem: see~\cite{KN13} for an overview. Smooth Fano polytopes -- corresponding to toric Fano manifolds -- can be classified efficiently using an algorithm of {\O}bro~\cite{Obr07}. The number grows slowly with dimension: just~$18$ out of the~$105$ three-dimensional Fano manifolds are toric, and even by dimension eight there are only~$749\,892$ smooth Fano polytopes.~\emph{Reflexive} Fano polytopes -- corresponding to Gorenstein toric Fano varieties -- are of particular importance in mirror symmetry: topologically mirror-symmetric pairs of Calabi--Yau varieties can be constructed as hypersurfaces in Gorenstein toric Fano varieties, and this duality is reflected in the combinatorics of the polytopes~$P$ and~$P^*$~\cite{Bat94}. There are~$4319$ reflexive polytopes in dimension three, and~$473\,800\,776$ in dimension four~\cite{KS98,KS00}. Toric Fano varieties with only terminal singularities correspond to Fano polytopes~$P\subset\NQ$ such that the only lattice points in~$P$ are the origin and the vertices. These have been classified in dimension three~\cite{Kas06}: there are~$634$ cases. Toric Fano varieties with canonical singularities correspond to Fano polytopes~$P\subset\NQ$ such that the only interior lattice point of~$P$ is the origin. These~\emph{canonical Fano polytopes} have also been classified in dimension three~\cite{Kas10,Zenodo-canonical3}: there are~$674\,688$ cases. 

These polytope classifications provide a natural place to search for rigid MMLPs. For example, an analysis of the~$4319$ reflexive three-dimensional polytopes was used in~\cite{ACGK12,CCGK16} to recover the~$98$ deformation families of three-dimensional Fano manifolds with very ample anticanonical bundle. Extending this idea, by constructing rigid MMLPs supported on Fano polytopes more generally, we can begin to populate the geography of \QQFano{} threefolds. This is developed further in~\cite{CHKP22,Heu22}, but let us first consider how much of the geography in Figure~\ref{fig:grdb} we might hope to realise this way.

\begin{figure}[tbp]
\centering
\begin{subfigure}{.48\textwidth}
\includegraphics[width=\linewidth]{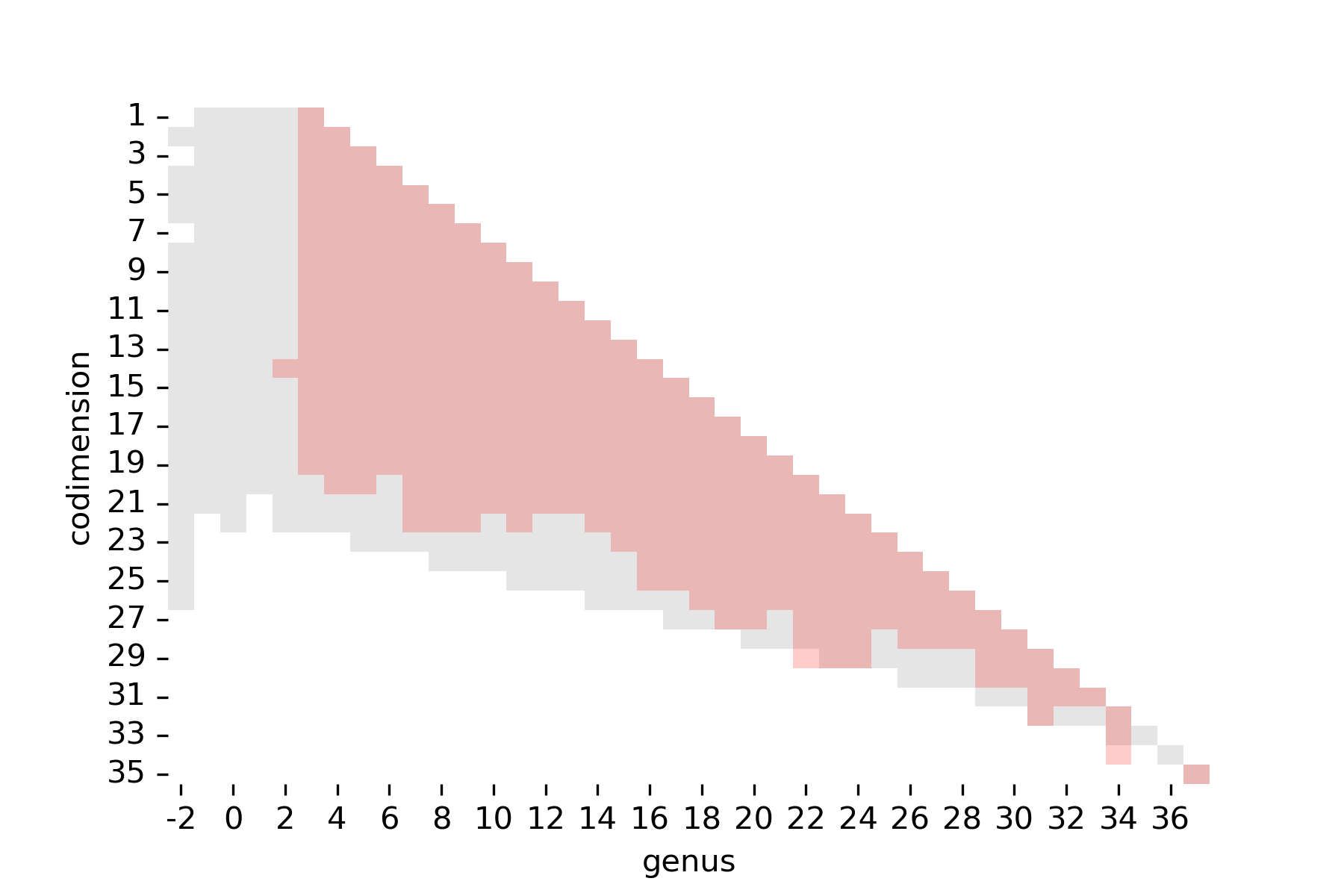}
\caption{Canonical polytopes ($1$ interior point).}
\label{fig:coverage_1pt}
\end{subfigure}
\hfill
\begin{subfigure}{.48\textwidth}
\includegraphics[width=\linewidth]{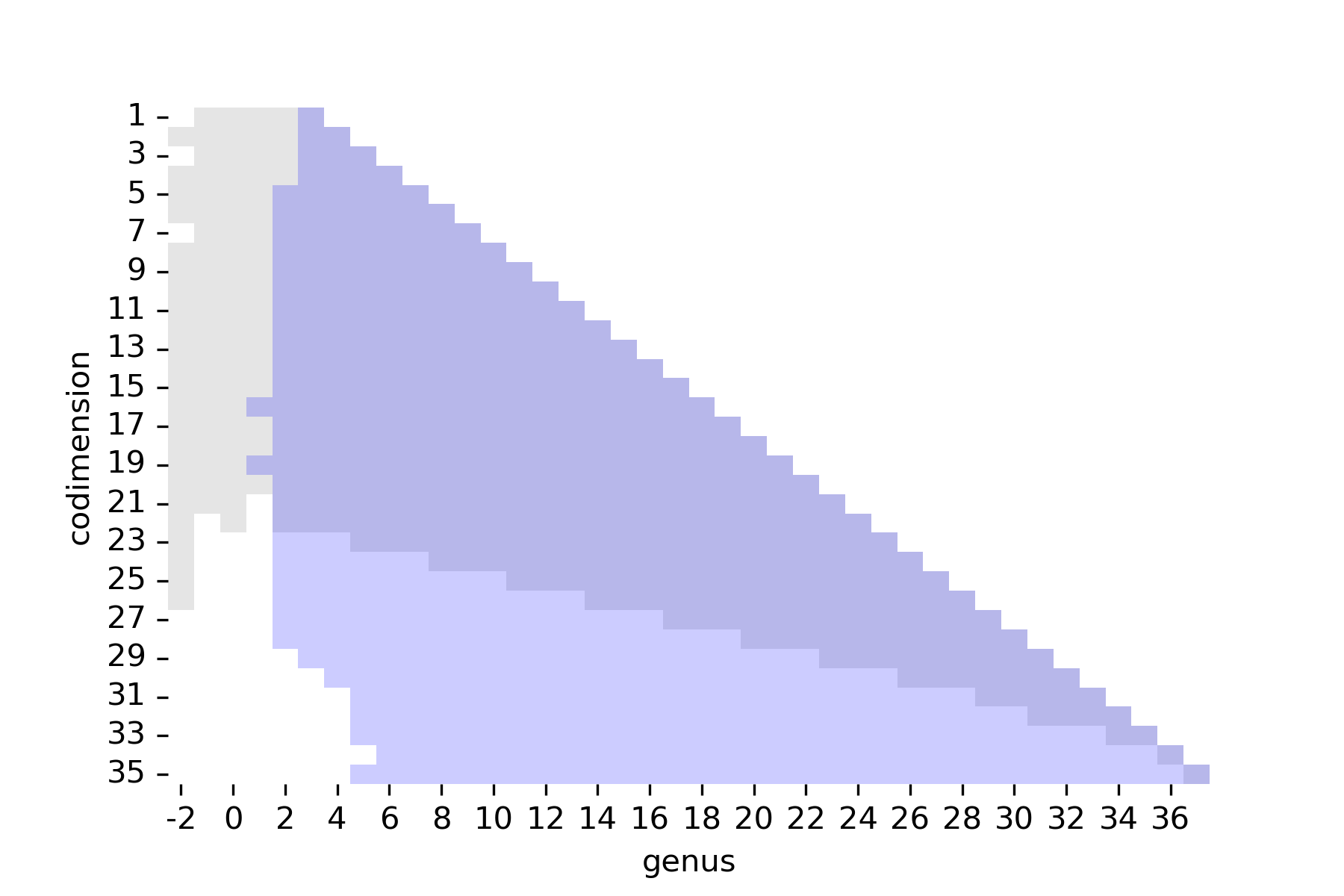}
\caption{Fano polytopes with~$2$ interior points.}
\label{fig:coverage_2pt}
\end{subfigure} 
\caption{The portions of the landscape of \QQFano{} threefolds that might be realised by three-dimensional Fano polytopes with one and two interior lattice points, overlaid on Figure~\ref{fig:grdb}.}
\label{fig:coverage}
\end{figure}

\begin{figure}[tbp]
\includegraphics[width=0.7\textwidth]{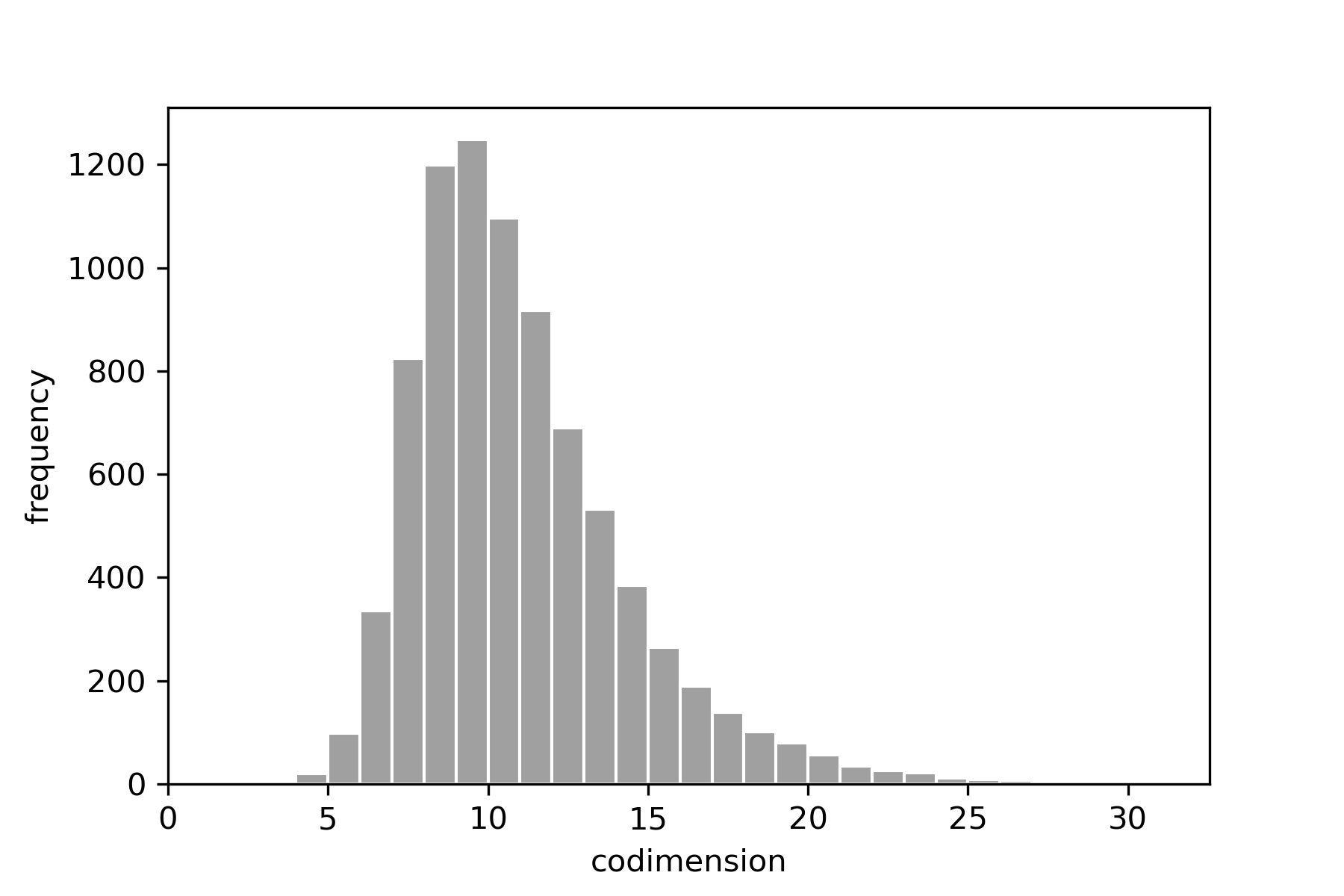}
\caption{Known mutation-equivalence classes of rigid MMLPs supported on $3$-dimensional canonical Fano polytopes, by codimension.}
\label{fig:codims}
\end{figure}

Genus and codimension are deformation invariants. We expect, as discussed above, that a \QQFano{} threefold~$X$ corresponding to a rigid MMLP~$f$ admits a degeneration to the toric variety~$X_f$ defined by the spanning fan of the Newton polytope $P\subset\NQ$ of~$f$. In this case we can read off the genus and codimension of~$X_f$ (and hence of~$X$) from the combinatorics of the dual polytope $P^*$ of $P$. Namely, the genus is given by $|P^*\cap M|-2$, and the codimension is computed from the Hilbert basis of the four-dimensional cone $\RR_{\geq 0}\cdot(P^*\times\{1\})$. Figure~\ref{fig:coverage_1pt} shows the pairs~$(g(Y), \codim{Y})$ realised by toric Fano varieties~$Y$ with at worst canonical singularities\footnote{See~\cite[Figure~6]{BK22} for a more detailed analysis.}; here~$Y$ is playing the role of~$X_f$. It should be emphasised that not every canonical Fano polytope will support a rigid MMLP: initial experiments indicate~${\sim}15\%$ of the canonical Fano polytopes support at least one rigid MMLP, resulting in approximately~$8300$ mutation-equivalence classes of rigid MMLPs. The codimensions of the hypothetical \QQFano{} threefolds discovered (or rather, suggested) in this way are summarised in Figure~\ref{fig:codims}: note that the majority of these have quite high codimension.

We expect that increasing the number of interior lattice points allowed in the Newton polytope (equivalently, allowing log-terminal singularities in the corresponding toric Fano varieties) will reveal new mutation-equivalence classes of rigid MMLPs. For example, Fano polytopes with two interior lattice points are classified in~\cite{BK16}. The pairs~$(g(Y), \codim{Y})$ for the resulting toric Fano varieties~$Y$ are shown in Figure~\ref{fig:coverage_2pt}. Once again we expect only a fraction of these Fano polytopes to support a rigid MMLP, and hence degenerate to a \QQFano{} threefold. In particular, we expect that none of the two-point polytopes that fall outside the gray shaded region in Figure~\ref{fig:coverage_2pt} -- or equivalently the region shown in Figure~\ref{fig:grdb} -- support a rigid MMLP. Figure~\ref{fig:coverage} suggests that increasing the number of interior lattice points might allow us to find \QQFano{} threefolds in both higher codimension and lower genus. The latter point is intuitively plausible: as the Newton polytope~$P$ grows larger, so the dual polytope~$P^*$ grows smaller and so we expect~$|P^* \cap M|$ and hence the corresponding genus to decrease.
%-------------------------------------------------------------------------------
\section{The $95$ hypersurfaces}\label{sec:95}
%-------------------------------------------------------------------------------
We return to the~$95$ \QQFano{} threefold hypersurfaces $X_d$ in weighted projective space
discussed in~\S\ref{sec:low_codim}.  As in~\S\ref{sec:mirror_symmetry}, we can apply the Givental/Hori--Vafa construction to obtain a Laurent polynomial mirror for~$X_d$. Write

\begin{align*}
 W = x_0 + x_1 + x_2 + x_3 + x_4 && \text{where} &&&
 x_0 x_1^{a_1} x_2^{a_2} x_3^{a_3} x_4^{a_4} = 1 \\
 && &&& x_1 + x_2 + x_3 + x_4 = 1.
\end{align*}
Eliminating the second equation by setting 
\begin{align*}
 x_1 = \frac{1}{1+x+y+z} &&
 x_2 = \frac{x}{1+x+y+z} &&
 x_3 = \frac{y}{1+x+y+z} &&
 x_4 = \frac{z}{1+x+y+z}
\end{align*}
and solving for~$x_0$ yields the Laurent polynomial
\begin{align*}
 f = W - 1 - c = \frac{(1+x+y+z)^d}{x^{a_2}y^{a_3}z^{a_4}} - c&&\text{where}&&c=\frac{d!}{a_1! a_2! a_3! a_4!}
\end{align*}
Here the constant~$c$ is chosen to ensure that the coefficient~$c'_1$ of the classical period~\eqref{eq:pi} is zero. The corresponding Newton polytope~$P$ is a simplex; the number of interior lattice points is given by the tetrahedral number~$(d-3)(d-2)(d-1)/6$. 

The hypersurface with the largest value for~$d$ is~$X_{66}\subset\PP(1,5,6,22,33)$. It has genus~$-1$, and gives a simplex with~$43\,680$ interior lattice points. Experiments mutating the Laurent polynomial show that we cannot expect to substantially decrease the number of interior lattice points. This agrees with our intuition in~\S\ref{sec:high_codim} that constructing low-genus \QQFano{} threefolds requires Fano polytopes with large numbers of interior points. Classifying Fano polytopes with tens of interior lattice points, let alone tens of thousands, is well beyond the reach of current methods, so it is implausible that we could find examples such as the 95 hypersurfaces using mirror symmetry methods alone.
%-------------------------------------------------------------------------------
\section{An example in codimension two}\label{sec:codim_2}
%-------------------------------------------------------------------------------
As discussed, we expect that if a \QQFano{} threefold~$X$ corresponds under mirror symmetry to a Laurent polynomial $f$, then $X$ admits a degeneration to the toric variety~$X_f$ defined by the spanning fan of the Newton polytope $P\subset\NQ$ of~$f$. If this is the case, then only those \QQFano{} threefolds that admit toric degenerations can be detected using mirror symmetry methods. For example, consider the codimension~two \QQFano{} threefold~$X_{12,14}\subset\PP(2,3,4,5,6,7)$, which is part of the Iano-Fletcher classification~\cite{Fle00}. This has~$\dim{H^0(X,{-K}_X)}=0$, and so can never deform to a toric variety\footnote{It would be interesting to know if in general the vanishing of~$H^0(X,{-K}_X)$ is the only obstruction to~$X$ admitting a toric degeneration.}. So mirror symmetry methods alone will not be enough to complete the classification.
%-------------------------------------------------------------------------------
\section{A two-pronged approach}\label{sec:2pronged}
%-------------------------------------------------------------------------------
Classical construction methods (\S\ref{sec:low_codim}, \S\ref{sec:95}--\ref{sec:codim_2}) and the mirror symmetry approach (\S\ref{sec:mirror_symmetry}--\ref{sec:high_codim}) complement each other well, with the classical methods most effective for \QQFano{} threefolds of low codimension and mirror symmetry giving remarkable coverage in high codimension. Mirror symmetry also gives new insight at one of the most difficult points of the analysis: determining whether a given set of possible numerical invariants is realised by a \QQFano{} variety. Conjecturally this problem is equivalent to determining the existence (or non-existence) of rigid MMLPs with appropriate Newton polytopes: a challenging but tractable combinatorial problem, which is amenable to machine computation. Traditionally there are no systematic methods available here.

This suggests a new approach -- a hybrid of classical and mirror symmetry methods -- to populate and analyse the landscape of \QQFano{} threefolds. First, use graded ring methods~\cite{grdb,ABR02,BK22} to determine an explicit list of possible Hilbert series. Then, for each candidate Hilbert series $H(t)$ on that list, classify rigid MMLPs $f$ supported on Fano polytopes~$P$ with dual Ehrhart series $\Ehr(P^*,t)$ equal to $H(t)$, where $P$ lies in an appropriate combinatorial class. This predicts the number of deformation families, as well as specific $\QQ$\nobreakdash-Gorenstein toric degenerations $X_f$ of these families. Finally, use deformation, Laurent inversion, or unprojection methods to construct each family. Meaningful analysis beyond particular cases will require cluster-scale computation and database methods, because of the numbers of cases involved (thousands of deformation families, hundreds of thousands of polytopes), but this is well within what is now possible. For a substantial step in this direction, see~\cite{CHKP22,Heu22,CKP22}.

We close by discussing a classical example from this new point of view. Consider a sequence of \QQFano{} threefolds~$Y^{(i)}\subset\PP(1,1,1,1,1,2,\ldots,2)$ in codimension $i+1$ that arise by successively unprojecting~$i$ planes inside a quartic hypersurface~$X_4\subset\PP^3$. The case when $i=1$ is discussed in~\S\ref{sec:unproj}. When~$i\ge3$, there may be more than one configuration of planes: see~\cite{BD22}. By determining mutation-equivalence classes of rigid MMLPs, we predict numbers of distinct deformation families as shown in Table~\ref{tab:num_families}. There remains much detailed analysis still to do, but at a coarse level this prediction matches an observation of Brown--Ducat~\cite{BD22}: that a quartic threefold can contain at most seven planes meeting only pairwise in points. Indeed a quartic $X$ that contains a plane $\Pi$ is birational to a cubic del~Pezzo fibration, obtained as the projection away from $\Pi$. Any other plane in $X$ that meets $\Pi$ only in points corresponds to a section of $(-1)$-curves in this fibration, and thus there are at most six such planes, corresponding to a maximal set of six disjoint lines in the general fibre of this cubic del~Pezzo fibration; cf.~\cite{Kal07}. Moreover a set of seven such planes can be explicitly realised inside the Burkhardt quartic~$X_4^{\textrm{Bur}}\subset\PP^4$, which is the unique quartic having the maximal number of~$45$ ordinary double points~\cite{DSV} and contains a configuration of~$40$ planes. Therefore, if we constrain the unprojection plane configurations to meet as simply as possible, there is no configuration of eight planes to unproject into codimension nine, in accordance with the conjecture.

\begin{table}[tbp]
    \centering
    \small
    \[
    \begin{array}{rccccccccc}
    \toprule
    i & 0 & 1 & 2 & 3 & 4 & 5 & 6 & 7 & 8 \\
    \midrule
    \text{GRDB ID} & \grdb{20521} & \grdb{20522} & \grdb{20543} & \grdb{20652} & \grdb{20938} & \grdb{21436} & \grdb{22093} & \grdb{22801} & \grdb{23386} \\
    \text{Codimension} & 1 & 2 & 3 & 4 & 5 & 6 & 7 & 8 & 9 \\
    \text{Number of families} & 1 & 1 & 1 & 3 & 7 & 11 & 11 & 4 & 0\\
    \bottomrule
    \end{array}
    \]
    \caption{The predicted number of deformation families of $\QQ$-Fano threefolds $Y^{(i)}$ from~\S\ref{sec:2pronged}.}
    \label{tab:num_families}
\end{table}
%-------------------------------------------------------------------------------
\bibliographystyle{plain}

\begin{thebibliography}{10}

\bibitem{ACCHKOPPT16}
Mohammad Akhtar, Tom Coates, Alessio Corti, Liana Heuberger, Alexander~M.
  Kasprzyk, Alessandro Oneto, Andrea Petracci, Thomas Prince, and Ketil
  Tveiten.
\newblock Mirror symmetry and the classification of orbifold del {P}ezzo
  surfaces.
\newblock {\em Proc. Amer. Math. Soc.}, 144(2):513--527, 2016.

\bibitem{ACGK12}
Mohammad Akhtar, Tom Coates, Sergey Galkin, and Alexander~M. Kasprzyk.
\newblock Minkowski polynomials and mutations.
\newblock {\em SIGMA Symmetry Integrability Geom. Methods Appl.}, 8:Paper 094,
  17, 2012.

\bibitem{Ale94}
Valery Alexeev.
\newblock General elephants of {${\bf Q}$}-{F}ano 3-folds.
\newblock {\em Compositio Math.}, 91(1):91--116, 1994.

\bibitem{Alt98}
Selma Alt{\i}nok.
\newblock {\em Graded Rings Corresponding to Polarised {K}3 Surfaces and
  $\mathbb{Q}$-{F}ano 3-folds}.
\newblock PhD thesis, University of Warwick, 1998.

\bibitem{ABR02}
Selma Alt{\i}nok, Gavin Brown, and Miles Reid.
\newblock Fano 3-folds, {$K3$} surfaces and graded rings.
\newblock In {\em Topology and geometry: commemorating {SISTAG}}, volume 314 of
  {\em Contemp. Math.}, pages 25--53. Amer. Math. Soc., Providence, RI, 2002.

\bibitem{AGM94}
Paul~S. Aspinwall, Brian~R. Greene, and David~R. Morrison.
\newblock Calabi-{Y}au moduli space, mirror manifolds and spacetime topology
  change in string theory.
\newblock {\em Nuclear Phys. B}, 416(2):414--480, 1994.

\bibitem{BK16}
Gabriele Balletti and Alexander~M. Kasprzyk.
\newblock Three-dimensional lattice polytopes with two interior lattice points.
\newblock \href{https://arxiv.org/abs/1612.08918}{arXiv:1612.08918 [math.CO]},
  2016.

\bibitem{Bat94}
Victor~V. Batyrev.
\newblock Dual polyhedra and mirror symmetry for {C}alabi-{Y}au hypersurfaces
  in toric varieties.
\newblock {\em J. Algebraic Geom.}, 3(3):493--535, 1994.

\bibitem{Bel}
Pieter Belmans.
\newblock Fanography: a tool to visually study the geography of {F}ano 3-folds.
\newblock \url{https://www.fanography.info}.

\bibitem{Bir21}
Caucher Birkar.
\newblock Singularities of linear systems and boundedness of {F}ano varieties.
\newblock {\em Ann. of Math. (2)}, 193(2):347--405, 2021.

\bibitem{BCHM10}
Caucher Birkar, Paolo Cascini, Christopher~D. Hacon, and James McKernan.
\newblock Existence of minimal models for varieties of log general type.
\newblock {\em J. Amer. Math. Soc.}, 23(2):405--468, 2010.

\bibitem{BD22}
Gavin Brown and Tom Ducat.
\newblock Double unprojection of {F}ano 3-folds.
\newblock In preparation.

\bibitem{grdb}
Gavin Brown and Alexander~M. Kasprzyk.
\newblock The {G}raded {R}ing {D}atabase.
\newblock \url{http://www.grdb.co.uk}.

\bibitem{BK16b}
Gavin Brown and Alexander~M. Kasprzyk.
\newblock Four-dimensional projective orbifold hypersurfaces.
\newblock {\em Exp. Math.}, 25(2):176--193, 2016.

\bibitem{fanodata}
Gavin Brown and Alexander~M. Kasprzyk.
\newblock The {F}ano 3-fold database.
\newblock \emph{Zenodo} \url{https://doi.org/10.5281/zenodo.5820338}, 2022.

\bibitem{BK22}
Gavin Brown and Alexander~M. Kasprzyk.
\newblock Kawamata boundedness for {F}ano threefolds and the {G}raded {R}ing
  {D}atabase.
\newblock \href{https://arxiv.org/abs/2201.07178}{arXiv:2201.07178 [math.AG]},
  2022.

\bibitem{BKQ18}
Gavin Brown, Alexander~M. Kasprzyk, and Muhammad~Imran Qureshi.
\newblock Fano 3-folds in {$\Bbb P^2\times \Bbb P^2$} format, {T}om and
  {J}erry.
\newblock {\em Eur. J. Math.}, 4(1):51--72, 2018.

\bibitem{BKZ19}
Gavin Brown, Alexander~M. Kasprzyk, and Lei Zhu.
\newblock Gorenstein formats, canonical and {C}alabi--{Y}au threefolds.
\newblock {\em Experimental Mathematics}, pages 1--19, 2019.

\bibitem{BKR12}
Gavin Brown, Michael Kerber, and Miles Reid.
\newblock Fano 3-folds in codimension 4, {T}om and {J}erry. {P}art {I}.
\newblock {\em Compos. Math.}, 148(4):1171--1194, 2012.

\bibitem{CF93}
F.~Campana and H.~Flenner.
\newblock Projective threefolds containing a smooth rational surface with ample
  normal bundle.
\newblock {\em J. Reine Angew. Math.}, 440:77--98, 1993.

\bibitem{COGP91}
Philip Candelas, Xenia~C. de~la Ossa, Paul~S. Green, and Linda Parkes.
\newblock A pair of {C}alabi-{Y}au manifolds as an exactly soluble
  superconformal theory.
\newblock {\em Nuclear Phys. B}, 359(1):21--74, 1991.

\bibitem{CP17}
Ivan Cheltsov and Jihun Park.
\newblock Birationally rigid {F}ano threefold hypersurfaces.
\newblock {\em Mem. Amer. Math. Soc.}, 246(1167):v+117, 2017.

\bibitem{CCC11}
Jheng-Jie Chen, Jungkai~A. Chen, and Meng Chen.
\newblock On quasismooth weighted complete intersections.
\newblock {\em J. Algebraic Geom.}, 20(2):239--262, 2011.

\bibitem{CCGGK13}
Tom Coates, Alessio Corti, Sergey Galkin, Vasily Golyshev, and Alexander~M.
  Kasprzyk.
\newblock Mirror symmetry and {F}ano manifolds.
\newblock In {\em European {C}ongress of {M}athematics}, pages 285--300. Eur.
  Math. Soc., Z\"{u}rich, 2013.

\bibitem{CCGK16}
Tom Coates, Alessio Corti, Sergey Galkin, and Alexander~M. Kasprzyk.
\newblock Quantum periods for 3-dimensional {F}ano manifolds.
\newblock {\em Geom. Topol.}, 20(1):103--256, 2016.

\bibitem{CHKP22}
Tom Coates, Liana Heuberger, and Alexander~M. Kasprzyk.
\newblock {M}irror symmetry, {L}aurent inversion and the classification of
  $\mathbb{Q}$-{F}ano threefolds.
\newblock \href{https://arxiv.org/abs/2210.07328}{arXiv:2210.07328 [math.AG]},
  2022.

\bibitem{CKP22}
Tom Coates, Alexander~M. Kasprzyk, and Giuseppe Pitton.
\newblock Computing maximally mutable {L}aurent polynomials.
\newblock In preparation.

\bibitem{CKPT21}
Tom Coates, Alexander~M. Kasprzyk, Giuseppe Pitton, and Ketil Tveiten.
\newblock Maximally mutable {L}aurent polynomials.
\newblock {\em Proc. of the Royal Society A.}, 477(2254):Paper No. 20210584,
  21, 2021.

\bibitem{CKP15}
Tom Coates, Alexander~M. Kasprzyk, and Thomas Prince.
\newblock Four-dimensional {F}ano toric complete intersections.
\newblock {\em Proc. of the Royal Society A.}, 471(2175):Paper No. 20140704,
  14, 2015.

\bibitem{CKP19}
Tom Coates, Alexander~M. Kasprzyk, and Thomas Prince.
\newblock Laurent inversion.
\newblock {\em Pure Appl. Math. Q.}, 15(4):1135--1179, 2019.

\bibitem{Cor95}
Alessio Corti.
\newblock Factoring birational maps of threefolds after {S}arkisov.
\newblock {\em J. Algebraic Geom.}, 4(2):223--254, 1995.

\bibitem{CPR00}
Alessio Corti, Aleksandr Pukhlikov, and Miles Reid.
\newblock Fano {$3$}-fold hypersurfaces.
\newblock In {\em Explicit birational geometry of 3-folds}, volume 281 of {\em
  London Math. Soc. Lecture Note Ser.}, pages 175--258. Cambridge Univ. Press,
  Cambridge, 2000.

\bibitem{CR02}
Alessio Corti and Miles Reid.
\newblock Weighted {G}rassmannians.
\newblock In {\em Algebraic geometry}, pages 141--163. de Gruyter, Berlin,
  2002.

\bibitem{CD20}
Stephen Coughlan and Tom Ducat.
\newblock Constructing {F}ano 3-folds from cluster varieties of rank 2.
\newblock {\em Compos. Math.}, 156(9):1873--1914, 2020.

\bibitem{Dan78}
V.~I. Danilov.
\newblock The geometry of toric varieties.
\newblock {\em Uspekhi Mat. Nauk}, 33(2(200)):85--134, 247, 1978.

\bibitem{DSV}
A.~J. de~Jong, N.~I. Shepherd-Barron, and A.~Van~de Ven.
\newblock On the {B}urkhardt quartic.
\newblock {\em Math. Ann.}, 286(1-3):309--328, 1990.

\bibitem{delPez87}
Pasquale Del~Pezzo.
\newblock Sulle superficie dell'$n^{\text{mo}}$ ordine immerse nello spazio ad
  $n$ dimensioni.
\newblock {\em Rend. del Circolo Mat. di Palermo}, 1:241--255, 1887.

\bibitem{Dol82}
Igor Dolgachev.
\newblock Weighted projective varieties.
\newblock In {\em Group actions and vector fields ({V}ancouver, {B}.{C}.,
  1981)}, volume 956 of {\em Lecture Notes in Math.}, pages 34--71. Springer,
  Berlin, 1982.

\bibitem{DH16}
Charles~F. Doran and Andrew Harder.
\newblock Toric degenerations and {L}aurent polynomials related to {G}ivental's
  {L}andau-{G}inzburg models.
\newblock {\em Canad. J. Math.}, 68(4):784--815, 2016.

\bibitem{Duc18}
Tom Ducat.
\newblock Constructing {$\Bbb Q$}-{F}ano 3-folds \`a la {P}rokhorov \& {R}eid.
\newblock {\em Bull. Lond. Math. Soc.}, 50(3):420--434, 2018.

\bibitem{Fan47}
Gino Fano.
\newblock Nuove ricerche sulle variet\`a algebriche a tre dimensioni a
  curve-sezioni canoniche.
\newblock {\em Pont. Acad. Sci. Comment.}, 11:635--720, 1947.

\bibitem{Ful93}
William Fulton.
\newblock {\em Introduction to toric varieties}, volume 131 of {\em Annals of
  Mathematics Studies}.
\newblock Princeton University Press, Princeton, NJ, 1993.
\newblock The William H. Roever Lectures in Geometry.

\bibitem{Giv98}
Alexander Givental.
\newblock A mirror theorem for toric complete intersections.
\newblock In {\em Topological field theory, primitive forms and related topics
  ({K}yoto, 1996)}, volume 160 of {\em Progr. Math.}, pages 141--175.
  Birkh{\"a}user Boston, Boston, MA, 1998.

\bibitem{GP90}
B.~R. Greene and M.~R. Plesser.
\newblock Duality in {C}alabi-{Y}au moduli space.
\newblock {\em Nuclear Phys. B}, 338(1):15--37, 1990.

\bibitem{HM10}
Christopher~D. Hacon and James McKernan.
\newblock Existence of minimal models for varieties of log general type. {II}.
\newblock {\em J. Amer. Math. Soc.}, 23(2):469--490, 2010.

\bibitem{Heu22}
Liana Heuberger.
\newblock $\mathbb{Q}$-{F}ano threefolds and {L}aurent inversion.
\newblock \href{https://arxiv.org/abs/2202.04184}{arXiv:2202.04184 [math.AG]},
  2022.

\bibitem{HV00}
Kentaro Hori and Cumrun Vafa.
\newblock Mirror symmetry.
\newblock \href{https://arxiv.org/abs/hep-th/0002222}{arXiv:hep-th/0002222},
  2000.

\bibitem{Fle00}
A.~R. Iano-Fletcher.
\newblock Working with weighted complete intersections.
\newblock In {\em Explicit birational geometry of 3-folds}, volume 281 of {\em
  London Math. Soc. Lecture Note Ser.}, pages 101--173. Cambridge Univ. Press,
  Cambridge, 2000.

\bibitem{Isk77}
V.~A. Iskovskih.
\newblock Fano threefolds. {I}.
\newblock {\em Izv. Akad. Nauk SSSR Ser. Mat.}, 41(3):516--562, 717, 1977.

\bibitem{Isk78}
V.~A. Iskovskih.
\newblock Fano threefolds. {II}.
\newblock {\em Izv. Akad. Nauk SSSR Ser. Mat.}, 42(3):506--549, 1978.

\bibitem{Isk79}
V.~A. Iskovskih.
\newblock Anticanonical models of three-dimensional algebraic varieties.
\newblock In {\em Current problems in mathematics, {V}ol. 12 ({R}ussian)},
  pages 59--157, 239 (loose errata). VINITI, Moscow, 1979.

\bibitem{IM71}
V.~A. Iskovskih and Ju.~I. Manin.
\newblock Three-dimensional quartics and counterexamples to the {L}\"{u}roth
  problem.
\newblock {\em Mat. Sb. (N.S.)}, 86(128):140--166, 1971.

\bibitem{JK01}
Jennifer~M. Johnson and J{\'a}nos Koll{\'a}r.
\newblock Fano hypersurfaces in weighted projective 4-spaces.
\newblock {\em Experiment. Math.}, 10(1):151--158, 2001.

\bibitem{Kal07}
Anne-Sophie Kaloghiros.
\newblock {\em The topology of terminal quartic 3-folds}.
\newblock PhD thesis, University of Cambridge, 2007.

\bibitem{Kas06}
Alexander~M. Kasprzyk.
\newblock Toric {F}ano three-folds with terminal singularities.
\newblock {\em Tohoku Math. J. (2)}, 58(1):101--121, 2006.

\bibitem{Kas10}
Alexander~M. Kasprzyk.
\newblock Canonical toric {F}ano threefolds.
\newblock {\em Canad. J. Math.}, 62(6):1293--1309, 2010.

\bibitem{Zenodo-canonical3}
Alexander~M. Kasprzyk.
\newblock The classification of toric canonical {F}ano 3-folds.
\newblock \emph{Zenodo} \url{https://doi.org/10.5281/zenodo.5866330}, 2010.

\bibitem{KN13}
Alexander~M. Kasprzyk and Benjamin Nill.
\newblock Fano polytopes.
\newblock In {\em Strings, gauge fields, and the geometry behind}, pages
  349--364. World Sci. Publ., Hackensack, NJ, 2013.

\bibitem{KNP17}
Alexander~M. Kasprzyk, Benjamin Nill, and Thomas Prince.
\newblock Minimality and mutation-equivalence of polygons.
\newblock {\em Forum Math. Sigma}, 5:Paper No. e18, 48, 2017.

\bibitem{KOW20}
In-Kyun Kim, Takuzo Okada, and Joonyeong Won.
\newblock K-stability of birationally superrigid {F}ano 3-fold weighted
  hypersurfaces.
\newblock \href{https://arxiv.org/abs/2011.07512}{arXiv:2011.07512 [math.AG]},
  2020.

\bibitem{KMM92}
J\'{a}nos Koll\'{a}r, Yoichi Miyaoka, and Shigefumi Mori.
\newblock Rational connectedness and boundedness of {F}ano manifolds.
\newblock {\em J. Differential Geom.}, 36(3):765--779, 1992.

\bibitem{KMMT00}
J\'{a}nos Koll\'{a}r, Yoichi Miyaoka, Shigefumi Mori, and Hiromichi Takagi.
\newblock Boundedness of canonical {$\bold Q$}-{F}ano 3-folds.
\newblock {\em Proc. Japan Acad. Ser. A Math. Sci.}, 76(5):73--77, 2000.

\bibitem{KS98}
Maximilian Kreuzer and Harald Skarke.
\newblock Classification of reflexive polyhedra in three dimensions.
\newblock {\em Adv. Theor. Math. Phys.}, 2(4):853--871, 1998.

\bibitem{KS00}
Maximilian Kreuzer and Harald Skarke.
\newblock Complete classification of reflexive polyhedra in four dimensions.
\newblock {\em Adv. Theor. Math. Phys.}, 4(6):1209--1230, 2000.

\bibitem{Mor75}
Shigefumi Mori.
\newblock On a generalization of complete intersections.
\newblock {\em J. Math. Kyoto Univ.}, 15(3):619--646, 1975.

\bibitem{Mor82}
Shigefumi Mori.
\newblock Threefolds whose canonical bundles are not numerically effective.
\newblock {\em Ann. of Math. (2)}, 116(1):133--176, 1982.

\bibitem{Mor88}
Shigefumi Mori.
\newblock Flip theorem and the existence of minimal models for {$3$}-folds.
\newblock {\em J. Amer. Math. Soc.}, 1(1):117--253, 1988.

\bibitem{MM81}
Shigefumi Mori and Shigeru Mukai.
\newblock Classification of {F}ano {$3$}-folds with {$B\sb{2}\geq 2$}.
\newblock {\em Manuscripta Math.}, 36(2):147--162, 1981/82.

\bibitem{MM03}
Shigefumi Mori and Shigeru Mukai.
\newblock Erratum: ``{C}lassification of {F}ano 3-folds with {$B_2\geq 2$}''.
\newblock {\em Manuscripta Math.}, 110(3):407, 2003.

\bibitem{Obr07}
Mikkel {\O}bro.
\newblock An algorithm for the classification of smooth {F}ano polytopes.
\newblock \href{https://arxiv.org/abs/0704.0049}{arXiv:0704.0049 [math.CO]},
  2007.

\bibitem{PR04}
Stavros~Argyrios Papadakis and Miles Reid.
\newblock Kustin-{M}iller unprojection without complexes.
\newblock {\em J. Algebraic Geom.}, 13(3):563--577, 2004.

\bibitem{Pol05}
Joseph Polchinski.
\newblock {\em String theory. {V}ol. {II}}.
\newblock Cambridge Monographs on Mathematical Physics. Cambridge University
  Press, Cambridge, 2005.
\newblock Superstring theory and beyond, reprint of 2003 edition.

\bibitem{Pri18}
Thomas Prince.
\newblock Smoothing toric {F}ano surfaces using the {G}ross-{S}iebert
  algorithm.
\newblock {\em Proc. Lond. Math. Soc. (3)}, 117(3):617--660, 2018.

\bibitem{Pro10}
Yuri Prokhorov.
\newblock {$\Bbb Q$}-{F}ano threefolds of large {F}ano index, {I}.
\newblock {\em Doc. Math.}, 15:843--872, 2010.

\bibitem{PR16}
Yuri Prokhorov and Miles Reid.
\newblock On {$\Bbb Q$}-{F}ano 3-folds of {F}ano index 2.
\newblock In {\em Minimal models and extremal rays ({K}yoto, 2011)}, volume~70
  of {\em Adv. Stud. Pure Math.}, pages 397--420. Math. Soc. Japan, [Tokyo],
  2016.

\bibitem{Prz11}
Victor Przyjalkowski.
\newblock Hori-{V}afa mirror models for complete intersections in weighted
  projective spaces and weak {L}andau-{G}inzburg models.
\newblock {\em Cent. Eur. J. Math.}, 9(5):972--977, 2011.

\bibitem{QS11}
Muhammad~Imran Qureshi and Bal\'{a}zs Szendr\H{o}i.
\newblock Constructing projective varieties in weighted flag varieties.
\newblock {\em Bull. Lond. Math. Soc.}, 43(4):786--798, 2011.

\bibitem{Rei80}
Miles Reid.
\newblock Canonical {$3$}-folds.
\newblock In {\em Journ\'{e}es de {G}\'{e}ometrie {A}lg\'{e}brique d'{A}ngers,
  {J}uillet 1979/{A}lgebraic {G}eometry, {A}ngers, 1979}, pages 273--310.
  Sijthoff \& Noordhoff, Alphen aan den Rijn---Germantown, Md., 1980.

\bibitem{Rei02}
Miles Reid.
\newblock Graded rings and birational geometry.
\newblock In {\em Proceedings of algebraic geometry symposium ({K}inosaki,
  {O}ct 2000), K. Ohno (Ed.)}, pages 1--72. 2002.
\newblock Available from \url{http://www.warwick.ac.uk/~masda/3folds/}.

\bibitem{San95}
Takeshi Sano.
\newblock On classifications of non-{G}orenstein {${\bf Q}$}-{F}ano {$3$}-folds
  of {F}ano index {$1$}.
\newblock {\em J. Math. Soc. Japan}, 47(2):369--380, 1995.

\bibitem{San96}
Takeshi Sano.
\newblock Classification of non-{G}orenstein {${\bf Q}$}-{F}ano {$d$}-folds of
  {F}ano index greater than {$d-2$}.
\newblock {\em Nagoya Math. J.}, 142:133--143, 1996.

\bibitem{Tak02}
Hiromichi Takagi.
\newblock On classification of {$\Bbb Q$}-{F}ano 3-folds of {G}orenstein index
  2. {I}, {II}.
\newblock {\em Nagoya Math. J.}, 167:117--155, 157--216, 2002.

\bibitem{Tak89}
Kiyohiko Takeuchi.
\newblock Some birational maps of {F}ano {$3$}-folds.
\newblock {\em Compositio Math.}, 71(3):265--283, 1989.

\end{thebibliography}

%-------------------------------------------------------------------------------
\end{document}